\documentclass[twoside,11pt]{amsart} \usepackage{epsfig}
\usepackage{latexsym} 
\usepackage{amsmath} \usepackage{amssymb}



\newtheorem{theorem}{Theorem}[section]

\newtheorem{conjecture}[theorem]{Conjecture}

\newdimen\epsfxsize
\newdimen\epsfysize

\newcommand {\gap}     {\makebox[0.075 in]{}}   
\newcommand {\st}      {\gap : \gap}   

\newcommand {\set}[1]  {\left\{ {#1} \right\}}   
\newcommand {\ord}[1]  {{#1}^{\rm th}}

\newcommand{\primeprod}[1] {{#1}^\#}

\newcommand{\Z}     {{\mathbb Z}}

\newcommand{\N}[2]  {N_{{#2}}({#1})}

\newcommand{\Est}[2]  {E_{{#2}}({#1})}
\newcommand{\Cnt}[2] {C_{{#2}}({#1})}

\newcommand{\pgap}   {{\mathcal G}}

\begin{document}

\title{Estimating constellations among primes \\ I - Uniformity}

\date{20 Mar 2008; 4 Dec 2013}

\author{Fred B. Holt and Helgi Rudd}
\address{fbholt@u.washington.edu ;  4311-11th Ave NE \#500, Seattle, WA 98105;
Unit 5 / 270 Campbell Parade, Bondi Beach, 2026, Sydney, Australia}

\begin{abstract}
A few years ago we identified a recursion that works directly with the gaps among
the generators in each stage of Eratosthenes sieve.  This recursion provides explicit
enumerations of sequences of gaps among the generators, which are known as 
constellations.

In this paper, we use those enumerations to estimate the numbers of these constellations
that occur as constellations among prime numbers, and we compare these estimates
with computational results.  We include in our estimates the constellations corresponding
to three and four consecutive primes in arithmetic progression.

For these initial estimates, we assume that the copies of a given constellation tend toward 
a uniform distribution in the cycle of gaps, as the recursion progresses.
Our simple estimates based on the recursion of gaps and the assumption of uniformity 
appear to have correct asymptotic
behavior, and they exhibit a systematic error correlated to length of the constellation.
\end{abstract}

\maketitle

\section{Introduction}
We work with the prime numbers in ascending order, denoting the
$\ord{k}$ prime by $p_k$.  Accompanying the sequence of primes
is the sequence of gaps between consecutive primes.
We denote the gap between $p_k$ and $p_{k+1}$ by
$g_k=p_{k+1}-p_k.$
These sequences begin
$$
\begin{array}{rrrrrrc}
p_1=2, & p_2=3, & p_3=5, & p_4=7, & p_5=11, & p_6=13, & \ldots\\
g_1=1, & g_2=2, & g_3=2, & g_4=4, & g_5=2, & g_6=4, & \ldots
\end{array}
$$

A number $d$ is the {\em difference} between prime numbers if there are
two prime numbers, $p$ and $q$, such that $q-p=d$.  There are already
many interesting results and open questions about differences between
prime numbers; a seminal and inspirational work about differences
between primes is Hardy and Littlewood's 1923 paper \cite{HL}.

A number $g$ is a {\em gap} between prime numbers if it is the difference
between consecutive primes; that is, $p=p_i$ and $q=p_{i+1}$ and
$q-p=g$.
Differences of length $2$ or $4$ are also gaps; so open questions
like the Twin Prime Conjecture, that there are an infinite number
of gaps $g_k=2$, can be formulated as questions about differences
as well.

A {\em constellation among primes} \cite{Riesel} is a sequence of consecutive gaps
between prime numbers.  Let $s=a_1 a_2 \cdots a_k$ be a sequence of $k$
numbers.  Then $s$ is a constellation among primes if there exists a sequence of
$k+1$ consecutive prime numbers $p_i p_{i+1} \cdots p_{i+k}$ such
that for each $j=1,\ldots,k$, we have the gap $p_{i+j}-p_{i+j-1}=a_j$.  
Equivalently,
$s$ is a constellation if for some $i$ and all $j=1,\ldots,k$,
$a_j=g_{i+j}$.

We will write the constellations without marking a separation between
single-digit gaps.  For example, a constellation of $24$ denotes
a gap of $g_k=2$ followed immediately by a gap $g_{k+1}=4$.
For the small primes we will consider
explicitly, most of these gaps are single digits, and the separators introduce
a lot of visual clutter.  We use commas only to separate double-digit gaps in
the cycle.  For example, a constellation of $2,10,2$ denotes a gap of $2$
followed by a gap of $10$, followed by another gap of $2$.

In \cite{FBHgaps} Holt introduced a recursion that works directly on the gaps
among the generators in each stage of Eratosthenes sieve.
These are the generators of $Z \bmod \primeprod{p}$ in which $\primeprod{p}$
is the product of the prime numbers up to $p$.
For a constellation $s$, this recursion enables us to enumerate exactly how
many copies of $s$ occur in the $\ord{k}$ stage of the sieve.  We denote this
number of copies of $s$ as $\N{p_k}{s}$.

For example, after the primes $2$, $3$, and $5$ and their multiples have been
removed, we have the cycle of gaps $\pgap(\primeprod{5}) = 64242462$.  This
cycle of $8$ gaps sums to $30$.  In this cycle, for the constellation $s=2$, we 
have $\N{5}{2}=3$.  For the constellation $s=242$, we have $\N{5}{242}=1$.
The cycle of gaps $\pgap(\primeprod{p})$ has $\phi(\primeprod{p})$ gaps that sum
to $\primeprod{p}$.

In this paper we use the counts of constellations among the generators of
$Z \bmod \primeprod{p}$ to estimate the numbers of constellations that occur
among prime numbers.  For $Z \bmod \primeprod{p_k}$, all the generators that
lie between $p_{k+1}$ and $p_{k+1}^2$ are actually prime numbers, so the
corresponding constellations are actually constellations among prime numbers.
For our estimates we assume uniform distributions of the various constellations
in the cycle of gaps.

We make estimates for the single-gap constellations $s=2$ (twin primes), 
$s=6$, and $s=8$.
Our estimates exceed the true counts for these constellations
 by around $16\%$ up to $1 E+12$.
Over this range, one of the estimates by Hardy and Littlewood \cite{HL} for $s=2$ 
underestimates the true counts for twin primes by around $8\%$.
Regarding twin primes, we also make estimates for the constellations
$s=242$ (prime quadruplets \cite{NicelyQuads}), $s=2,10,2$, and $s=2,10,2,10,2$.
Finally, to explore arithmetic sequences of prime numbers, which have been visible
in the recent work of Green and Tao
 \cite{BGreen},  we make estimates for $s=66$ and $s=666$.
These correspond respectively to arithmetic progressions of three and four consecutive
primes.  

For the few constellations we have considered, our estimates appear
to have the correct asymptotic behavior, but our estimates also seem to have a
systematic error correlated with the number of gaps in the constellation.
For a constellation $s$ of $j$ gaps, our estimates of how often $s$ occurs as
a constellation among prime numbers exceed the actual counts by approximately
$9\% \times (j+1)$ when our computations reach $1.0E+12$.

The recursion itself is interesting, since it provides exact counts \cite{FBHgaps} of constellations
among the generators $Z \bmod \primeprod{p}$, as a representation of Eratosthenes sieve.
In this paper, we use simple first-order statistics to make our estimates, and over the
range of our computations, these estimates have the right asymptotic behavior.
We believe that more sophisticated statistical models for the distributions of 
constellations in the cycle of gaps will produce more accurate estimates.

\section{Related Results}
There are of course several avenues of research into the distribution of
primes.  Much research into constellations has been motivated by
two conjectures:  the twin primes conjecture, and Hardy and Littlewood's
broader $k$-tuple conjecture \cite{HL,Rib}.

The twin primes conjecture asserts that the gap $g=2$ occurs infinitely often.
Work on this conjecture
has included computer-based enumerations \cite{IJ,NicelyTwins,PSZ}
and investigations of Brun's constant \cite{Sieves,Riesel,HL,Rib}.
Brun's constant is the sum of the reciprocals of twin primes.
This series is known to converge, 
and the sharpest current estimate \cite{Rib} is $1.902160577783278.$
One generalization of the twin primes conjecture is a conjecture by Polignac from
1849 \cite{Rib} that for every even positive integer $N$ 
there are an infinite number of gaps $g_k=N$.

Most estimates for sequences of differences, e.g. \cite{HL,MV}, are derived by treating probabilities on residues as independent probabilities.  
Computational confirmation of these estimates has been carried out
by several researchers, notably in \cite{Brent3,NicelyTwins}.
Some researchers have applied their investigations of differences among
primes to study gaps and constellations.  The general surveys \cite{Rib,Riesel}
provide overviews of some of this work, and the estimate
from the Hardy-Littlewood paper \cite{HL} can be used for constellations consisting of
$2$'s and $4$'s, e.g. $24$, $42$, $242$, $424$, etc.
Brent \cite{Brent} applied the principle of inclusion and exclusion
to those estimates to obtain strong estimates for the
gaps $2,4,6,\ldots,80$.
Richards \cite{Rich} conjectured
that the constellation $24$ occurs infinitely often.
Clement \cite{quads} and Nicely \cite{NicelyQuads} have addressed the
constellation $242$, corresponding to {\em prime quadruplets}, pairs of twin
primes separated by a gap of $4$.
 
In contrast, the recursion identified in \cite{FBHgaps}
 preserves the structure in the cycles of gaps at each stage of Eratosthenes sieve.  The occurrence of constellations in stages of the sieve is entirely deterministic.  We
 make estimates for complicated constellations based on assuming uniform distributions
 of the constellations for which we can obtain exact enumerations via the recursion.

Our results are not quite commensurate with the tables and examples previously
published.
The tables and examples of \cite{Brent,Brent3,IJ,NicelyTwins} provide
estimates and counts of gaps with respect to large powers of ten.
Since we work directly with Eratosthenes sieve, our estimates
are given with respect to intervals $[p,p^2]$ for primes $p$.

\section{Recursion on Cycle of Gaps}
In the cycle of gaps, the first gap corresponds to the next prime.  In $\pgap(\primeprod{5})$
the first gap $g_1=6$, which is the gap between $1$ and the next prime, $7$.  The next
several gaps are actually gaps between prime numbers.  In the cycle of gaps
$\pgap(\primeprod{p_k})$, 
the gaps between $p_{k+1}$ and $p_{k+1}^2$ are in fact gaps between prime numbers.

There is a simple recursion which generates $\pgap(\primeprod{p_{k+1}})$ from
$\pgap(\primeprod{p_k})$.  This recursion and many of its properties are developed
in \cite{FBHgaps}.  We include only the concepts and results we need for the estimates in
this paper.

The recursion on the cycle of gaps consists of three steps.
\begin{itemize}
\item[R1.] The next prime $p_{k+1} = g_1+1$, one more than the first gap;
\item[R2.] Concatenate $p_{k+1}$ copies of $\pgap(\primeprod{p_k})$;
\item[R3.] Add adjacent gaps as indicated by the elementwise product 
$p_{k+1}*\pgap(\primeprod{p_k})$:  let $i_1=1$ and add together $g_{i_1}+g_{i_1+1}$; then for 
$n=1,\ldots,\phi(N)$, add $g_{j}+g_{j+1}$ and let 
$i_{n+1}=j$ if the running sum of the concatenated gaps from $g_{i_n}$ to $g_j$ is
$p_{k+1}*g_{n}.$
\end{itemize}

\noindent{\bf Example: $\pgap(\primeprod{7})$.}
To illustrate this recursion,
we construct $\pgap(\primeprod{7})$ from $\pgap(\primeprod{5})=64242462$.

\begin{itemize}
\item[R1.] Identify the next prime, $p_{k+1}= g_1+1 = 7.$
\item[R2.] Concatenate seven copies of $\pgap(\primeprod{5})$:
$$64242462 \; 64242462 \; 64242462 \; 64242462\; 64242462 \; 64242462 \;64242462$$
\item[R3.] Add together the gaps after the leading $6$ and 
thereafter after differences of 
$ 7*\pgap(\primeprod{5}) = 42, 28, 14, 28, 14, 28, 42, 14 $:
\begin{eqnarray*}
\pgap(\primeprod{7}) 
 &=&{\scriptstyle 
  6+\overbrace{\scriptstyle 424246264242}^{42}+
 \overbrace{\scriptstyle 4626424}^{28}+\overbrace{\scriptstyle 2462}^{14}+
 \overbrace{\scriptstyle 6424246}^{28}+\overbrace{\scriptstyle 2642}^{14}+
 \overbrace{\scriptstyle 4246264}^{28}+\overbrace{\scriptstyle 242462642424}^{42}+
 \overbrace{\scriptstyle 62 \; }^{14}} \\
 &=& {\scriptstyle 
 10, 2424626424 6 62642 6 46 8 42424 8 64 6 24626 6 4246264242, 10, 2}
\end{eqnarray*}
The final difference of $14$ wraps around the end of the cycle,
 from the addition preceding the final $6$ to the 
addition after the first $6$.
\end{itemize}

\subsection{Key properties of the recursion}
We summarize a few properties of the cycle of gaps $\pgap(\primeprod{p})$, as established in 
\cite{FBHgaps}.  The cycle of gaps ends in a $2$, and except for this final $2$, the cycle of gaps
is symmetric.  In constructing $\pgap(\primeprod{p_{k+1}})$,
each possible addition of adjacent gaps in $\pgap(\primeprod{p_k})$ occurs exactly once.

For a constellation $s$ of $j$ gaps in $\pgap(\primeprod{p_k})$, if $j < p_{k+1}-1$, then
in $\pgap(\primeprod{p_{k+1}})$ there will be at least
$p_{k+1}-j-1$ copies corresponding to this particular constellation $s$.
Since the additions in step R3 are spaced according to the
element-wise product $p_{k+1}*\pgap(\primeprod{p_k})$,
the minimal distance between additions is $2p_{k+1}$.
If the sum of the gaps in $s$, denoted $\sigma(s)$, is less than $2p_{k+1}$, then
the $j+1$ additions in step R3 that affect $s$ will occur in distinct copies of $s$
from step R2, and there will be exactly $p_{k+1}-j-1$ copies of $s$.

\subsection{Numbers of constellations}
The power of the recursion on the cycle of gaps is seen in the following theorem,
which enables us to calculate the number of occurrences of a constellation $s$ through
successive stages of Eratosthenes sieve.

\begin{theorem}\label{CountThm}
(from \cite{FBHgaps})
Let $s$ be a constellation of $j$ gaps in $\pgap(\primeprod{p_k})$, such that 
$j < p_{k+1}-1$ and  $\sigma(s) < 2p_{k+1}$.
Let $S$ be the set of all constellations $\bar{s}$ which would produce
$s$ upon one addition of differences.
Then the number $\N{p}{s}$ of occurrences of $s$ in $\pgap(\primeprod{p})$
 satisfies the recurrence
$$\N{p_{k+1}}{s} = (p_{k+1}-(j+1)) \cdot \N{p_k}{s}
 + \sum_{\bar{s} \in S} \N{p_k}{\bar{s}}$$
\end{theorem}

\begin{figure}[tb] 
\centering
\includegraphics[width=5in]{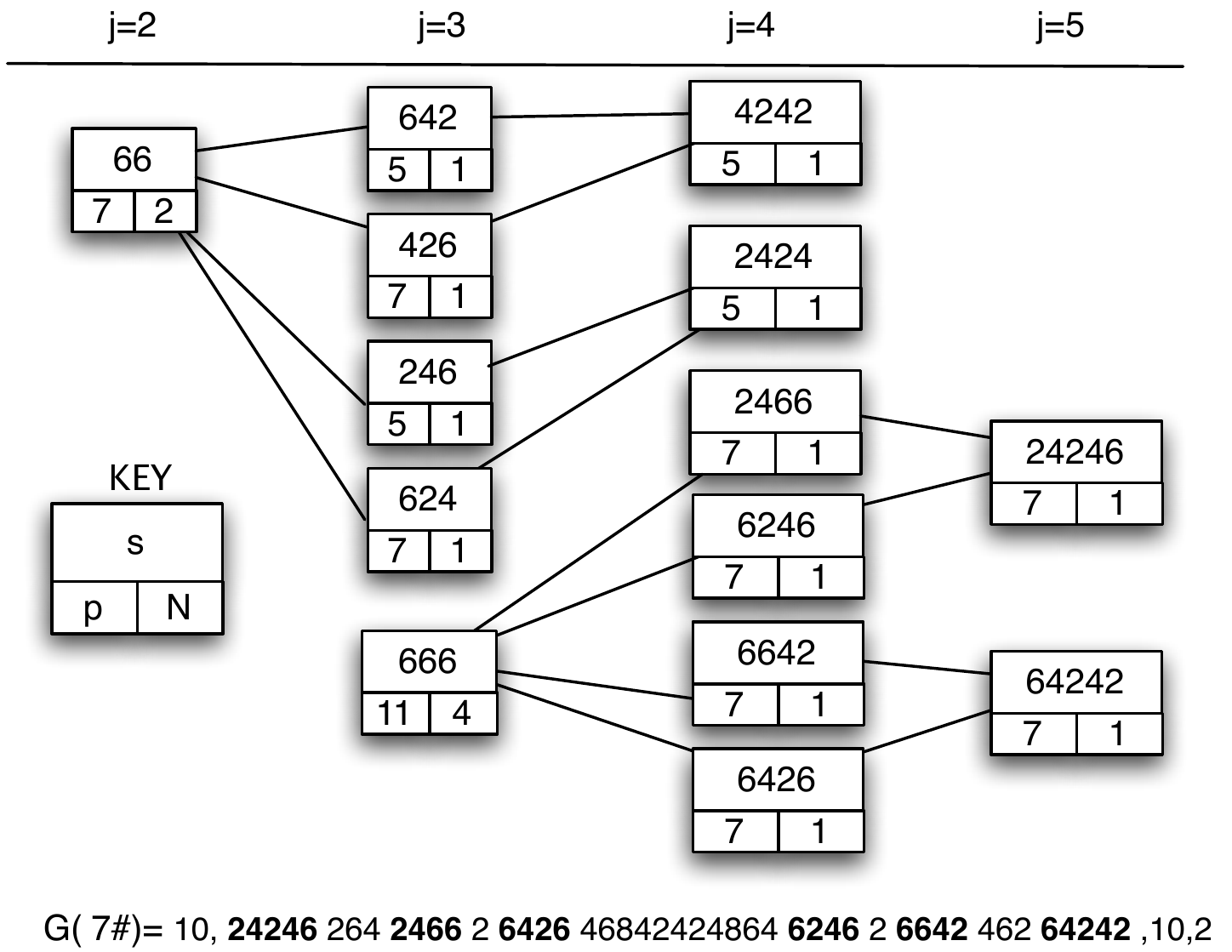}
\caption{\label{GrowthFig} This figure illustrates the initial conditions and driving terms for
calculating the numbers of copies of the constellations $s=66$ and $s=666$ in 
$\pgap(\primeprod{p})$.  These constellations correspond respectively
 to arithmetic sequences of three and four consecutive primes.  The entries in this chart
 indicate the 
 constellation $s$, its length $j$,
the prime for which the constellation occurs in $\pgap(\primeprod{p})$ and which satisfies the
conditions of Theorem~\ref{CountThm}, and the number $N=\N{p}{s}$ of occurrences of the
constellation in $\pgap(\primeprod{p})$.  From these figures we can derive the recursive count
$\N{q}{s}$, and from these counts we can make estimates $\Est{q}{s}$ 
for primes $q > p$.  At the bottom of the figure we show 
$\pgap(\primeprod{7})$ with the drivers of $s=666$ shown in bold.}
\end{figure}

Figure~\ref{GrowthFig} illustrates the initial conditions for the constellations $s=66$ and
$s=666$ and their driving terms.  Note that the initial conditions are not predicated on
when the constellations first appear but on the $\pgap(\primeprod{p})$ for which
the constellations satisfy the conditions of Theorem~\ref{CountThm}.

\section{Estimating constellations among prime numbers}
Based on the recursion above, we conjecture that the copies of a constellation $s$ 
are distributed approximately uniformly in the cycle of gaps $\pgap(\primeprod{p})$.
Under this assumption, we can estimate the number of
copies that occur in the interval $[p_{k},p_{k}^2]$:
\begin{equation}\label{EqEs}
\Est{k}{s} = \frac{p_{k}^2-p_{k}}{\primeprod{p_{k-1}}} \N{p_{k-1}}{s}.
\end{equation}

For a few interesting constellations, we compare these estimates to counts of how many
copies of $s$ actually occur in the interval $[q,q^2]$ for primes $q$.  We denote these
counts as $\Cnt{q}{s}$.

\subsection{Nearly Uniform Distributions}
The recursion in \cite{FBHgaps} suggests that if we track the images
of some constellation $s$ through several stages of Eratosthenes sieve, these images
will be almost uniformly distributed in the cycle of gaps.

In $\pgap(\primeprod{p_k})$ pick any constellation $s$ of length $j$, that satisfies the conditions
of Theorem~\ref{CountThm}.
Step R2 of the recursion creates $p_{k+1}$ copies of this constellation,
and these copies {\em are} distributed uniformly
in the interval $[1,\primeprod{p_{k+1}}]$.  Step R3 removes $j+1$ of these
copies of $s$ through a single addition in each of the removed copies.  

We don't know the distribution of these $j+1$ removed copies among the
$p_{k+1}$ uniformly distributed copies of $s$.  However, the removals must
preserve the symmetry of $\pgap(\primeprod{p})$, and as Eratosthenes sieve
proceeds, $j$ remains fixed as $p$ continues to grow, so the impact of the
$j+1$ removals on the uniform distribution among the $p$ copies diminishes.
As we continue applying the
recursion, the abundant images of $s$ are pushed toward uniformity by
step R2 and trimmed by step R3 away from uniformity but preserving symmetry.

These observations, about the effects of the recursion on the distribution
of occurrences of a particular constellation, indicate an approximately uniform distribution, 
but there is more work to be done in understanding the distribution of copies of $s$
in $\pgap(p)$.  In this paper we make estimates based on this assumption of uniformity.

\begin{conjecture}\label{uniconj}
Under the recursion on the cycle of gaps, all constellations $s$
in $\pgap(p)$  with $\sigma(s) < 2p$ tend toward
a uniform distribution in $\pgap(P)$ for all primes 
$P \gg p$.
\end{conjecture}

This conjecture suggests that the step R3 in the recursion has an approximately uniform effect on the uniform distribution of copies created by step R2.  That is, in the absence of any additional
information about the distribution of removals, each of the copies of $s$ has an equal chance 
$\frac{j+1}{p}$ of being removed.  If we focus on a fixed population of $N$ copies of $s$, then
we expect $(1-\frac{j+1}{p})N$ of these copies to survive step R3.

To bolster this conjecture, we can step back and look at the aggregate population of all constellations.  Not all constellations can accumulate in the middle of $\pgap(p)$.  
Some must fall
in $[p,p^2]$.  If the distribution for some particular constellation is forever biased strongly toward
the middle of the cycle, then the distributions of some other constellations must compensate for
this bias.  That is, if some constellations fall below the expected number of occurrences, then 
other constellations must exceed the expectations.

\section{Computational Results}
Our computations for counting constellations among prime numbers 
are currently performed using a modified segmented 
Eratosthenes sieve. Constellation counts for multiple constellations of varied length can be 
computed concurrently in a single pass through the Eratosthenes sieve. Constellations are 
defined as a two dimensional array, an additional one dimensional array of pointers is used to 
track the match position of constellations as prime gaps are identified by the Eratosthenes sieve. 
 
Constellation counts are stored in memory in a two dimensional array, one dimension for the 
counts and one dimension for each constellation being analyzed. 
When a constellation pattern is matched, counts for that constellation are incremented for all 
prime numbers from $\lceil \sqrt{p}\rceil$ to $p - \sigma$ where $p$ is the prime number currently 
identified by the Eratosthenes sieve and $\sigma$ is the sum of the gaps in the constellation. 
 When the sieve process is complete, all results held in the two dimensional array of constellation 
counts are streamed to a plain text file. 
 
\subsection{Limits for these computations}
Limits can be classified as memory limits, data type limits and execution time limits. 
 Because the computation is using a segmented sieve, the process of identifying primes does not 
require significant memory, the algorithm allows configuration of the segment size (which directly 
correlates to memory usage) with the intention of optimizing interaction with L2 cache memory. 
 
The current prime sieve code is limited to $9.0 E15$, which can easily be extended by altering data 
types. However if we are to extend these calculations, we need to employ a faster 
sieve. The current algorithm was chosen for its relative simplicity. 
Integration of the Sieve of Atkins, or an optimized version of an Eratosthenes Sieve is 
recommended.
 
The primary memory consideration is the two dimensional array which holds constellation 
counts. This array is currently of size $N \times n$ where $N$ is the number of integers being sieved 
and $n$ is the number of constellations being analyzed. The data type being used defines 
memory requirements and the choice of data type is determined by the constellation count totals 
expected. With a standard double datatype (8 bytes, allowing constellation count totals up to 
$9.0 E15$) and analysis of $n=10$ constellations concurrently to $1 E12$
 results in approximately 80Mb 
of memory usage. With the same scenario run up to $1 E15$, approximately 2.4 GB of memory is 
required. 
 
The current practical limit is execution time. Execution time is largely spent either sieving for 
primes, or incrementing constellation counts. The sieve time complexity is 
${\mathcal O}(N \log \log N)$. 
The time complexity for incrementing constellation counts is determined by the number of 
constellations being sieved and the number of matches for the constellations. 

\subsection{Comparison of estimates to counts}
As a first application of this approach, we estimate the number of
twin primes between $p$ and $p^2$.  From the recursion in Theorem~\ref{CountThm}
we derive a recursion for these estimates as well.
\begin{eqnarray*}
\Est{p_{k+1}}{2} &=& \frac{(p_{k+1}^2-p_{k+1}) }{\primeprod{p_k}} \N{p_k}{2} \\
&=& \frac{p_{k+1}^2-p_{k+1} }{p_k^2-p_k}\; \frac{p_k-2}{p_k}
 \;  \frac{p_k^2-p_k}{\primeprod{p_{k-1}}} \; \N{p_{k-1}}{2} \\
  & = & \frac{p_{k+1}^2-p_{k+1} }{p_k^2-p_k} \; \frac{p_k-2}{p_k} \; \Est{p_k}{2}.
\end{eqnarray*}

In the figures that follow, we plot the percentage error of our estimates compared to the actual
counts for various constellations.
$$ {\rm err} (s, p) = \frac{\Est{p}{s} - \Cnt{p}{s}}{\Cnt{p}{s}} $$
Here $\Cnt{p}{s}$ is the count of actual occurrences of the constellation $s$ in the interval
$\left[ p, p^2\right]$.

\begin{figure}[tb] 
\centering
\includegraphics[width=5in]{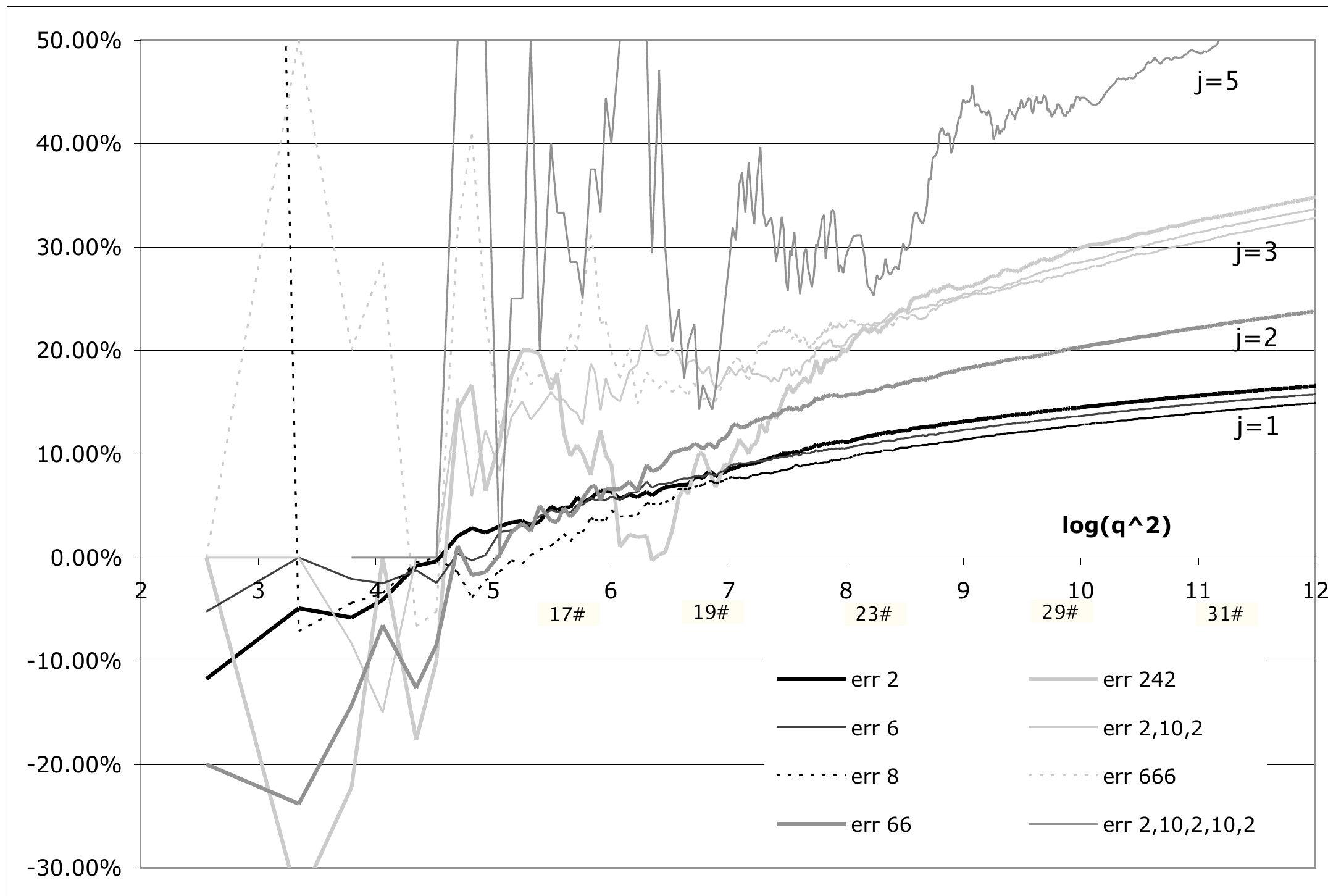}
\caption{\label{LogErrFig} This graph depicts the percentage error between our estimates and 
actual counts for various constellations.  Here we plot the percentage errors versus $\log q^2$. 
For the constellations we have considered, we observe that after initial noise,
our simple estimates have a systematic error that depends on the number $j$ of
gaps in the constellation.}
\end{figure}

Figure~\ref{LogErrFig} shows that after some initial noise, the errors settle down into trending curves.
We notice that the curves seem to sort out primarily by the length $j$ of the constellations.
The magnitude and duration of the noise seem also to depend on the length $j$ and the sum $\sigma$
of the gaps in the constellation.

The trends of the errors seem to depend primarily on the length $j$ of the constellation.
We observe from our computations that  the simple estimates appear to have the right
asymptotic behavior but that there are systematic errors that correlate to the number $j$
of gaps in the constellation.  
In Figure~\ref{LogErrFig} we note that our estimates for 
single-gap constellations are all off by around $16\%$.   
We have included
only one constellation with $j=2$ gaps, $s=66$, and the error through $1E+12$ is
around $24\%$.  
For the three constellations with $j=3$ gaps, $s=242$, $s=2,10,2$, and
$s=666$, our simple estimates all have an error of around $34\%$.  We have included
no constellations with $j=4$ symbols, and the single constellation with $j=5$, 
$s=2,10,2,10,2$ displays an error of $53\%$.  So there appears to be a systematic 
error of around $9\% \times (j+1)$ at $1E+12$, and these errors appear to be growing 
logarithmically. 
 We would like to find a model that accounts for this deviation between
the simple estimate and the actual counts.

\section{Observations on these estimates}
\subsection{Twin Primes} 
Our estimate for the number of twin primes must be contrasted with
those estimates provided by Hardy and Littlewood, which are supported by 
vast computation \cite{Brent,NicelyTwins}:
$$
\# \set{g_i = 2 \st p_i \in [2,N] }  \; \; \sim \;  2 c_2 \frac{N}{(\ln N)^2}
$$
in which the {\em twin prime constant} $c_2$ is given by the infinite product \cite{HL,Riesel} 
$c_2 = \prod \frac{p(p-2)}{(p-1)^2} = 0.6601618\ldots .$

\begin{figure}[tb] 
\centering
\includegraphics[width=5in]{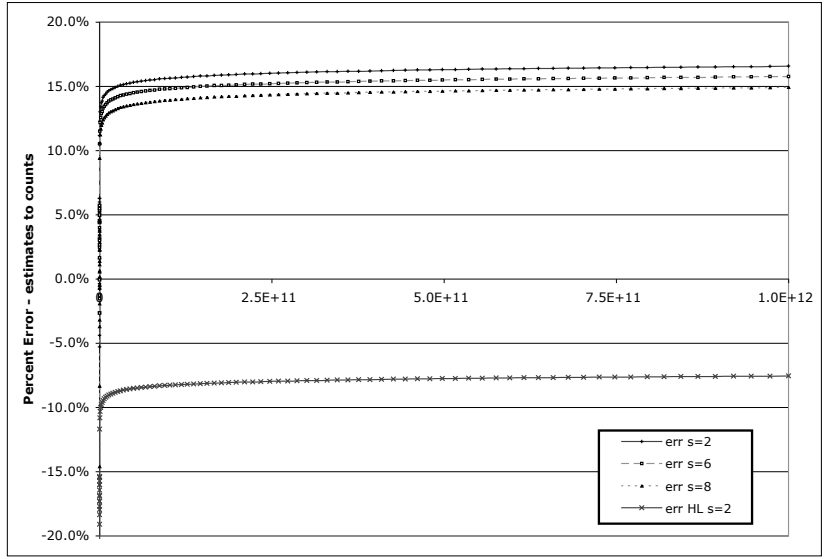}
\caption{\label{Pct268Fig} This chart depicts the percentage error between our estimates of the numbers of $2$'s, $6$'s, and $8$'s occurring as gaps between primes and actual counts for these gaps.  The counts and estimates are over the intervals $[q,q^2]$ for primes $q$.  
Here the percentage error is plotted against 
$q^2$, the top of the interval.  For $s=2$ we plot the percentage for both our estimate and
for the Hardy-Littlewood estimate ``HL 2".  On this linear scale, the early noise in the accuracy 
is all compressed against the vertical axis.}
\end{figure}

In Figure \ref{Pct268Fig} we plot the percentage errors between the estimates and 
the actual number of occurrences of the constellation $s=2$, for both
the Hardy-Littlewood estimate ``HL 2'' and our estimate $\Est{p}{2}$.
The comparisons are made over the intervals $[q,q^2]$ for primes $q$, and the percentage errors
are plotted against $q^2$.  Also depicted in Figure~\ref{Pct268Fig} are the percentage errors for
our estimates of the single-gap constellations $s=6$ and $s=8$.

Twin primes occur in interesting constellations, for example
the prime quadruplets constellation $s=242$.
This constellation occurs in $\pgap(5)$.
With $\N{5}{242}=1$, we calculate
the number of expected occurrences of the
constellation $242$ between $p$ and $p^2$.
Under the recursion, $\N{p_k}{242} = (p_k-4)\N{p_{k-1}}{242}$.
So, under our conjecture of uniformity, the expected number of these
constellations between $p_{k+1}$ and $p_{k+1}^2$ is
\begin{eqnarray*} 
\Est{p_{k+1}}{242} &=& \frac{p_{k+1}^2-p_{k+1}} {\primeprod{p_k}} \N{p_k}{242} \\
 & = & \frac{p_{k+1}^2-p_{k+1}}{p_k^2-p_k} \; \frac{p_k-4}{p_k} \;
  \frac{p_k^2-p_k} {\primeprod{p_{k-1}}} \N{p_{k-1}}{242} \\
 & = &  \frac{p_{k+1}^2-p_{k+1}}{p_k^2-p_k} \; \frac{p_k-4}{p_k} \Est{p_k}{242}.
 \end{eqnarray*}

With $s=242$, we can still compare to existing computations.
The extensive computations of \cite{NicelyQuads} support the Hardy-Littlewood
estimates \cite{Riesel} for $s=242$:
$$
\# \set{s = 242 \st s \subset [2,N]}  \;  \sim \; \frac{27}{2}c_4 \frac{N}{(\ln N)^4}$$
with $c_4 = \prod_{q\ge 5} \frac{q^3(q-4)}{(q-1)^4} = 0.30749\ldots .$
Figure~\ref{Pct666Fig} includes our estimate  $\Est{p}{242}$ and the Hardy-Littlewood
estimate ``HL 242''.

In contrast to approaches through differences, from the recursion on the cycle of gaps,
we can make estimates for 
the constellations $s=2,10,2$; $s=2,10,2,10,2$; $s=66$; and $s=666$.
Calculating $\Est{p}{2,10,2}$ involves driving
terms from the constellations $2462$ and $2642$.  We note that
for all $p\ge 13$, $\N{p}{2,10,2} > \N{p}{242}$, but that both have the
same dominant factor of $p-4$. 

\subsection{Primes in arithmetic progression}
Recently, Green and Tao \cite{BGreen} have offered a proof that there
exist arbitrarily long sequences of primes in arithmetic progression.
There are a number of researchers
who oversee enumerations and searches for primes and consecutive primes
in arithmetic progression.
In line with this avenue of research, we include in our computations and estimates the constellations
$s=66$ and $s=666$.   These constellations correspond to three and four 
consecutive primes in arithmetic progression.


A sequence of $k$ consecutive primes in arithmetic progression corresponds to
a constellation $s$ consisting of a single gap $g$ repeated $k-1$ times. 
For such a constellation $s$, the gap $g$ must be divisible by all primes $p \le k$.
For example, for $k=3$ the repeated gap must be divisible
by $6$, and $s=66$ is the minimal constellation for 
$k=3$ consecutive primes in arithmetic progression.
Similarly, $s=666$ is the minimal constellation for a $k=4$ 
consecutive primes in arithmetic progression, and for $k=5$ 
the minimal constellation is $s= 30,30,30,30$.

\begin{figure}[tb] 
\centering
\includegraphics[width=5in]{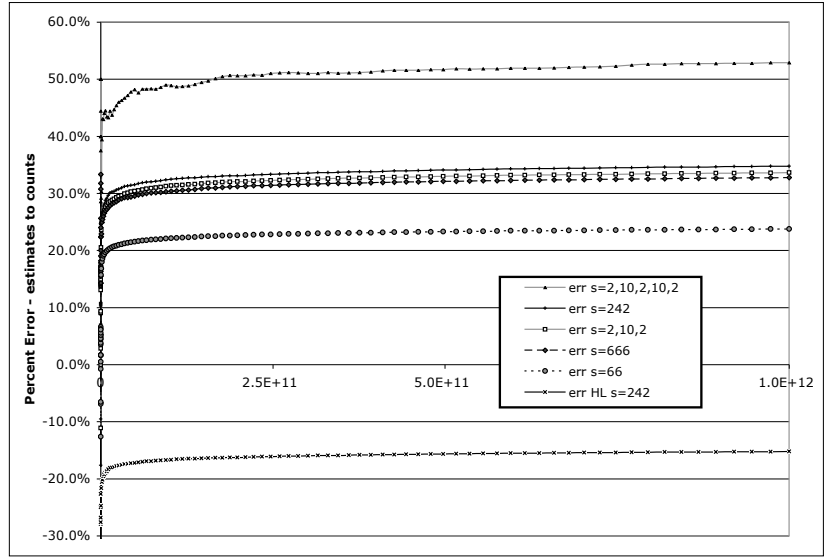}
\caption{\label{Pct666Fig} This chart depicts the percentage error between our estimates of 
the numbers of a few short constellations and the actual counts for these constellations.
The counts and estimates are over the intervals $[q,q^2]$ for primes $q$, and the
percentage error is plotted against $q^2$, the top of the interval.
The constellations $s=66$ and $s=666$ correspond to consecutive primes in 
arithmetic progression.
The constellation $s=242$ corresponds to prime quadruplets, 
and we plot the results of both our estimate and the Hardy-Littlewood estimate
``HL 242".  }
\end{figure}

We can apply Theorem~\ref{CountThm} and the other structure in the
recursion to these constellations.  
For this paper, we undertook enumerations and estimates for $s=66$ and $s=666$
up to $1E12$.
The constellation $s=66$ first occurs in $\pgap(\primeprod{7})$ and $s=666$ in 
$\pgap(\primeprod{11})$.   
The percentage errors for our simple estimates for $s=66$ and
$s=666$ are shown in Figure~\ref{Pct666Fig}.

Theorem~\ref{CountThm} provides conditions
under which a constellation corresponding to an arithmetic progression among possible 
primes flourishes within the cycles of gaps.
For  $k$ primes in an arithmetic progression,
if the corresponding constellation $s$ 
satisfies $\sigma(s) < 2p$ for any prime $p$ such that the arithmetic progression
falls in the interval $[p,\primeprod{p}]$, then the constellation $s$ survives
as a constellation among generators of $\Z \bmod q$ for all primes $q \ge p$.
Moreover, copies of $s$ undergo internal additions under the recursion,
and among the resulting constellations is $\bar{s}$ consisting of $k-1$
identical gaps.  By Theorem~\ref{CountThm} this constellation $\bar{s}$
persists in $\pgap(\primeprod{q})$ and its population grows superexponentially
with $q$.  

For example, both occurrences of 
$s=66$ in $\pgap(\primeprod{7})$ actually survive the recursion, corresponding to the
sequences of primes $47, 53, 59$ and $151, 157, 163$.
Using Theorem~\ref{CountThm} we know there will be
$26$ copies of $s=66$ in $\pgap(\primeprod{11})$.
Of these $26$ copies, the first four survive, 
as do three of the next nine copies.  Of the last $13$ copies
of $s=66$ in
$\pgap(\primeprod{11})$ none survive the recursion; each is subject
to at least one addition as the recursion continues.
In the interval $[13,2310]$ covered by $\pgap(\primeprod{11})$,
fifteen additional copies of $s=66$ are created from the driving terms of
$s=66$ (such as $\bar{s} = 624$).

The minimal constellation for $5$ consecutive primes in arithmetic progression
is $s=30,30,30,30$.
The sum $\sigma(s) = 120$, so the conditions of Theorem~\ref{CountThm}
are first satisfied for $\pgap(\primeprod{59})$.
The cycle $\pgap(\primeprod{59})$ consists of roughly $2.57E+20$
gaps of total sum around $1.92E+21$.  The constellation
$s=30,30,30,30$ first occurs among consecutive primes around ten million.

\subsection{Critiquing the simple estimate}
Based on our conjecture about uniformity, we have estimated the
number of copies of a constellation $s$ that fall between
$p_{k+1}$ and $p_{k+1}^2$ as
$$ \Est{p_{k+1}}{s} = (p_{k+1}^2-p_{k+1})\frac{\N{p_k}{s}}{\primeprod{p_k}}.$$

For small primes we can improve these estimates by adjusting the intervals for the
length of the constellation $s$ or the symmetry of the cycle of gaps $\pgap(\primeprod{p})$.
These adjustments change the estimates significantly for small primes, but their
impact fades rapidly in the face of factors of the form $p-j$.
From our computations, it seems that the estimates simply need to be adjusted by some
constant, which depends in part on the constellation $s$.  While we could estimate
the constant for the constellations we have enumerated, we do not yet have
a theoretical basis for adjusting the estimate in this way.

To make real improvement on the estimates, we need better statistical models for the
distribution of constellations in the cycle of gaps.  Since correlations among copies
of a constellation are largely preserved by the recursion, helpful second-order statistics
may be tractable.  Our conjecture~\ref{uniconj} does not take into account these correlations.
For example, the constellation $2,10,2,10,2$ occurs in $\pgap(7)$ and thereafter.
This constellation contains two occurrences of $2,10,2$ and
will consequently introduce jumps in $C_k(2,10,2)$. 

The structure of $\pgap(\primeprod{p})$ determines the constellations among primes
for intervals $[q, q^2]$ for primes $q$ up to $\sqrt{\primeprod{p}}$.   This means that any
deviation from the uniform distribution of constellations will persist at least through
$\sqrt{\primeprod{p}}$.

Even if our assumption of uniformity holds, the errors in our estimates may not be a surprise.
First of all, for our calculations out to $1E12$, we only need to repeat the complete recursion
out to $\pgap(\primeprod{37})$, and thereafter just perform the closures within this cycle of gaps
as the recursion continues.  The maximum prime $q$ with $q^2 < 37^\#$ is $2724079$.
So the calculations of our estimates use relatively small primes and may be biased by the 
relative abundance of small gaps early in the sequence of prime numbers.

The horizon of all current enumerations of twin primes and of prime quadruplets
\cite{Rib,PSZ,NicelyTwins, NicelyQuads} falls inside the cycle of gaps for
$p=53$.  $\pgap(\primeprod{53})$ has length $\primeprod{53} \approx 3.26 E19$.
That is, we have to run the complete recursion only up to  $\pgap(\primeprod{53})$,
and thereafter only perform the closures within $\pgap(\primeprod{53})$
up to $q=5708691479$, which is the maximum
prime $q$ with $q^2 < 53^\#$.

Second, we make our estimates from Equation (\ref{EqEs}) over the sample interval $[p_{k+1},p_{k+1}^2]$.
On the one hand, this interval is large enough for us to prove that the estimates for constellations of
lengths $1$ to $5$ gaps will grow with $p$, but the interval 
is a vanishingly small sample in $\pgap(\primeprod{p_k})$ as the process continues.  This interval
is also not a random sample from the distribution.  We are sampling from the start of the cycle of gaps, so
the errors will be correlated (e.g. $[p_k, p_k^2]$ and $[p_{k+1},p_{k+1}^2]$ substantially overlap).

Along the lines of this critique, although we are using $[p_{k+1},p_{k+1}^2]$ to 
sample $\pgap(\primeprod{p_k})$, the estimate would perhaps be better taken from
a cycle of gaps $\pgap(\primeprod{q})$ for a prime $q$ much smaller than $p_k$.
Perhaps we can reverse the framework for these estimates by fixing 
$\pgap(\primeprod{p_k})$ and estimating the survival of constellations within it as the
recursion continues.

\section{Conclusion}
We start with the recursion on the cycle of gaps in the stages of
Eratosthenes sieve.  Based on the recursion, we 
conjecture that all constellations, which occur in $\pgap(\primeprod{p})$
for some prime $p$ and the sum of whose gaps is less than $2p$,
tend toward a uniform distribution in later stages of the sieve.
From this conjecture, we estimate the number of
occurrences of a constellation between $p$ and $p^2$ for the
new prime $p$ at each stage of the sieve; all constellations
which occur before $p^2$ actually occur as constellations between primes.

We compare our estimates with enumerations for a sample of interesting constellations
with lengths from one gap to five gaps.  Our simple estimates compare well with the enumerations so far.
For single gaps or for constellations consisting only of $2$'s and $4$'s, 
other estimates are available \cite{HL, Brent, Riesel,NicelyQuads}, and 
our estimates agree with these established estimates to first order.

Having compared our simple estimates to actual counts through $1E+12$,
we observe a systematic error in the simple estimates, which depends on
the length of the constellation.  A focus of our continuing work in this area is
to understand this error and thereby to improve our estimates.


\bibliographystyle{amsplain}

\providecommand{\bysame}{\leavevmode\hbox to3em{\hrulefill}\thinspace}
\providecommand{\MR}{\relax\ifhmode\unskip\space\fi MR }
\providecommand{\MRhref}[2]{%
  \href{http://www.ams.org/mathscinet-getitem?mr=#1}{#2}
}
\providecommand{\href}[2]{#2}

\end{document}